\documentstyle{article}           
\pdfoutput=1 \pdfcompresslevel=9 \pdfadjustspacing=1
\font\eightrm = cmr10 scaled 800
\font\sixrm = cmr7 scaled 850
\font\fiverm = cmr5
\font\eightbf = cmbx10 scaled 800
\font\sixbf = cmbx7 scaled 850
\font\fivebf = cmbx5
\font\eightit = cmti10 scaled 800
\font\eighti = cmmi10 scaled 800
\font\sixi = cmmi7 scaled 850
\font\fivei = cmmi5
\font\eightsy = cmsy10 scaled 800
\font\sixsy = cmsy7 scaled 850
\font\fivesy = cmsy5
\font\eightsl = cmsl10 scaled 800
\font\eighttt = cmtt10 scaled 800
\font\fivefk = eufm5
\font\sixfk = eufm7 scaled 850
\font\sevenfk = eufm7
\font\eightfk = eufm10 scaled 800
\font\tenfk = eufm10
\newfam\fkfam
\textfont\fkfam=\tenfk \scriptfont\fkfam=\sevenfk
\scriptscriptfont\fkfam=\fivefk
\def\fk{\fam\fkfam}%
\font\fivebbm  = bbmsl10 scaled 500
\font\sixbbm   = bbmsl10 scaled 600
\font\sevenbbm = bbmsl10 scaled 700
\font\eightbbm = bbmsl10 scaled 800
\font\tenbbm   = bbmsl10
\newfam\bbmfam
\textfont\bbmfam=\tenbbm \scriptfont\bbmfam=\sevenbbm
\scriptscriptfont\bbmfam=\fivebbm
\def\bm{\fam\bbmfam}%
\def\eightpoint{%
\textfont0=\eightrm \scriptfont0=\sixrm \scriptscriptfont0=\fiverm
\def\rm{\fam0\eightrm}%
\textfont1=\eighti  \scriptfont1=\sixi  \scriptscriptfont1=\fivei
\textfont2=\eightsy \scriptfont2=\sixsy \scriptscriptfont2=\fivesy
\textfont\itfam=\eightit \def\it{\fam\itfam\eightit}%
\textfont\slfam=\eightsl \def\sl{\fam\slfam\eightsl}%
\textfont\ttfam=\eighttt \def\tt{\fam\ttfam\eighttt}%
\textfont\bffam=\eightbf \scriptfont\bffam=\sixbf
\scriptscriptfont\bffam=\fivebf \def\bf{\fam\bffam\eightbf}%
\textfont\fkfam=\eightfk \scriptfont\fkfam=\sixfk
\scriptscriptfont\fkfam=\fivefk
\textfont\bbmfam=\eightbbm \scriptfont\bbmfam=\sixbbm
\scriptscriptfont\bbmfam=\fivebbm
\rm}
\skewchar\eighti='177\skewchar\sixi='177
\skewchar\eightsy='60\skewchar\sixsy='60
\def\small{\eightpoint\baselineskip=9.6pt}%
\outer\def\title #1\by #2.{%
\pdfinfo{/Title (#1) /Author (#2) /Date (16 June 2014)}%
\centerline{{\useFnt[cmbx10 scaled 1440]#1}}%
\centerline{{\useFnt[cmr7](#2)}}%
\vskip 1ex plus .5ex}
\newcount\secNo \secNo=0
\outer\def\section #1.{%
\global\advance\secNo by 1
\global\ssecNo=0
\global\eqnNo=0
\goodbreak\vskip 3ex\noindent
{\useFnt[cmbx10 scaled 1200]\the\secNo.~#1}
\vglue 1ex}
\newcount\ssecNo \ssecNo=0
\outer\def\subsection #1.{%
\global\advance\ssecNo by 1
\goodbreak\vskip 2ex\noindent
{\useFnt[cmbx10]\the\secNo.\the\ssecNo.~#1}
\vglue 1ex}
\def\today{\number\day\space\ifcase\month\or
January\or February\or March\or April\or May\or June\or
July\or August\or September\or October\or November\or December\fi
\space\number\year}
\def\\{\hfill\break}
\catcode`@=11
\def\useFnt[#1]{\font\tmp@fnt = #1\relax\tmp@fnt}
\def\ifundef#1{\expandafter\ifx\csname#1\endcsname\relax}
\def\htag[#1]#2{\pdfdest name {#1} xyz\relax#2}%
\def\href[#1]#2{\leavevmode
\pdfstartlink
attr {/Border [0 0 0]} goto name {#1}\relax#2\pdfendlink}
\let\plaineqno\eqno
\newcount\eqnNo \eqnNo=0
\def\eqn@no{(\the\secNo.\the\eqnNo)}
\def\eqno#1$${\global\advance\eqnNo by 1
\plaineqno\htag[eqn.#1]{{\rm\eqn@no}}$$%
\ifundef{EQN@#1}%
\expandafter\xdef\csname EQN@#1\endcsname{\eqn@no}%
\else\errmessage{duplicated equation label}%
\fi\ignorespaces}
\def\reqn#1{\href[eqn.#1]{{\rm\csname EQN@#1\endcsname}}}
\newif\iflinenums  \linenumstrue 
\everydisplay={\textstyle}
\def\item@exec(#1)(#2,#3,#4)#5{{\parindent=#2\relax\par\leavevmode}%
\def\argA{#1}\relax\def\argB{#2}\relax%
\ifx\argA\argB\hangindent=#2\else\hangindent=#3\fi\hangafter=1%
\llap{#5\enspace}\ignorespaces}%
\def\item@optA[#1]#2{\item@exec(#1)(#1,,){#2}}%
\def\item@optB#1{\item@optA[20pt]{#1}}%
\def\item@opt{%
\ifx\item@char[\let\item@cmd\item@optA\relax\else%
\let\item@cmd\item@optB\fi\relax\item@cmd}%
\def\item{\futurelet\item@char\item@opt}

\newcount\foot@no \foot@no=0
\def\footnote{\global\advance\foot@no by 1
\edef\@sf{\spacefactor\the\spacefactor}%
\unskip{\raise.4\baselineskip\hbox{%
\useFnt[cmr7]\number\foot@no)}}\@sf
\insert\footins\bgroup\everydisplay={}\everypar={}\let\par\endgraf
\parindent=16pt \parskip=0pt \leftskip=0pt \rightskip=0pt
\splittopskip=10pt plus 1pt minus 1pt \floatingpenalty=20000
\ifundef{small}\else
\abovedisplayskip=6pt \abovedisplayshortskip=-4pt
\belowdisplayskip=4pt \belowdisplayshortskip=4pt
\small\fi
\smallskip
\textindent{\number\foot@no)}%
\bgroup\aftergroup\@foot\let\next}
\def\@foot{\egroup}
\skip\footins=12pt plus 2pt minus 4pt
\dimen\footins=.5\vsize
\newdimen\figwd \figwd=\hsize
\newcount\figno \figno=0
\def\put(#1,#2)#3{\unskip\leavevmode%
\vadjust{\smash{\raise#2\figwd\hbox{\hglue#1\figwd#3}}}
\ignorespaces}
\def\fig[#1]#2 -- #3 #4\endfig{\bgroup
\linenumsfalse\everypar={}\let\par\endgraf
\divide\figwd by 100\relax\figwd=#1\figwd
\vbox{\hsize=\figwd
\pdfximage width \figwd {#3}\pdfrefximage\pdflastximage
\relax #4
\vglue 0pt\centerline{#2}
}\egroup}
\catcode`@=12
\def\N{{\bm N}} \def\Z{{\bm Z}} \def\R{{\bm R}} \def\C{{\bm C}}
\def\RP{{\bm R}\!{\sl P}} \def\const{{\rm const}}
\let\em\it

\begin{document}
\title\vbox{\hsize=\hsize
 \hfil Discrete constant mean curvature nets in space forms:\hfil\vfil
 \hfil Steiner's formula and Christoffel duality\hfil}
\by A Bobenko${}^{\ast}$,%
 \insert\footins{\small\hglue 2em \llap{$\ast)$\enspace}%
  Partially supported by the DFG Collaborative Research
  Center TRR 109, ``Discretization in Geometry and Dynamics''.}
 U Hertrich-Jeromin \& I Lukyanenko.
\vskip 2em
\hbox to\hsize{\hfil\vbox{\hsize=0.8\hsize{\small{\bf Abstract.}\enspace
We show that the discrete principal nets in quadrics of constant
curvature that have constant mixed area mean curvature can be
characterized by the existence of a K\"onigs dual in a
concentric quadric.
\par}}\hfil}\vskip 1em
\hbox to\hsize{\hfil\vbox{\hsize=0.8\hsize{\small{\bf MSC 2010.}\enspace
{\it 53A10\/}, {\it 53C42\/}, 53A15, 53A30, 37K25, 37K35
\par}}\hfil}\vskip 1em
\hbox to\hsize{\hfil\vbox{\hsize=0.8\hsize{\small{\bf Keywords.}\enspace
K\"onigs net; K\"onigs dual; Christoffel transformation;
discrete surface; discrete conjugate net; space form;
hyperbolic space; projective space; M\"obius quadric;
constant mean curvature; mixed area.
\par}}\hfil}\vskip 1em

\section Introduction.
A variety of approaches have been pursued to obtain a notion of
discrete minimal and, more generally, discrete constant mean curvature
surfaces in Euclidean space and in spaces of constant curvature.  Two
fundamentally different starting points arise from the physical
interpretation of constant mean curvature surfaces as critical points
of an area functional and, on the other hand, from an integrable
systems viewpoint.  One principal difference of the two approaches is
that, in the former approach, the underlying combinatorial structure
is naturally that of a simplicial surface, cf [\href[ref.pipo93]{13,
Def 2}], whereas the integrable systems approach, discretizing
particular coordinate systems, demands a quadrilateral structure of
the discrete surface, cf [\href[ref.bopi99]{2}].

Restricting to quadrilateral meshes and, in particular, to planar
quadrilaterals (discrete conjugate nets) or quadrilaterals that have a
circumcircle (discrete curvature line nets), as we shall below, may at
first sight appear to be a severe restriction.  However, it turns out
that every surface patch can be approximated by discrete nets with
these properties, cf [\href[ref.bosu08]{4, Sects 5.3 and 5.6}],
reflecting the existence of the corresponding parametrizations in the
smooth case.  In fact, integrable discrete nets and discretizations of
their smooth counterparts in computer graphics are often visually
indistinguishable, reflecting a particular strength of these
integrable discretizations: the fact that these are {\it structure
preserving\/} discretizations, respecting key aspects of the geometry.
The apparent restriction of topology that is implied by the use of a
``global'' parameter net can also be overcome, for example, by using
more general ``quadgraphs'' ($2$-dimensional cell complexes with
quadrilateral $2$-cells) for a discrete domain manifold, where
singularities of the parameter net are modelled by vertices with
valences different from $4$, though this has so far only been
investigated in special cases, cf [\href[ref.bhs06]{3}].

Despite their differences --- regarding initial approach as well as
desired results --- the rich theory of constant mean curvature
surfaces in space forms suggests exploitation of central properties of
the surface class in either approach.  For example, the Lawson
correspondence between constant mean curvature surfaces in different
space forms has been crucial in the construction of simplicial
constant mean curvature surfaces in Euclidean space from simplicial
minimal surfaces in $S^3$, cf [\href[ref.gbpo96]{8, Sect 3}], as well
as in a definition of constant mean curvature nets in space forms via
their conformal deformation\footnote{The ``Lawson correspondence''
understood as a conformal deformation differs from the one between
constant mean curvature surfaces in $\R^3$ and minimal surfaces in
$S^3$ by a $90^\circ$-rotation that is commonly added to the latter.}
(Calapso transformation) as special isothermic nets, cf
[\href[ref.bjrs08]{6, Sect 5}].

In this paper we provide a simple geometric characterization of
discrete constant mean curvature surfaces in spaces of constant
curvature from an integrable systems point of view, as ``special''
discrete isothermic nets.

While there are no conceptual problems arising from the ambient
curvature $\kappa$ in the variational approach, see
[\href[ref.pipo93]{13, Def 3}], the original ideas to define discrete
constant mean curvature surfaces in Euclidean space from
[\href[ref.bopi96]{1, Sect 4}] and [\href[ref.jhp99]{9, Sect 5}]
cannot easily be generalized to curved ambient geometries.  A key
feature in the integrable systems approach to constant mean curvature
surfaces and to isothermic surfaces is the existence of B\"acklund and
Darboux transformations, respectively, alongside Bianchi permutability
theorems, which allow to construct discrete analogues of the smooth
surfaces by repeated transformations.  Characterizing (or defining)
constant mean curvature surfaces in terms of the isothermic
transformation theory, via the existence of distinguished transforms
as in [\href[ref.bopi96]{1}] and [\href[ref.jhp99]{9}], links in well
with the theory and allows to easily specialize the transformation
theory of discrete isothermic surfaces to discrete constant mean
curvature surfaces.  Using the conformal deformation of isothermic
nets, these characterizations of minimal and constant mean curvature
surfaces in Euclidean space lead to characterizations in space forms
via spherical or antipodal Darboux transforms --- as long as the
Lawson invariant $H^2+\kappa\geq0$: from [\href[ref.bjrs08]{6}] we
learn that these distinguished transforms are determined by solving a
quadratic equation with discriminant $H^2+\kappa$, so that in the case
$H^2+\kappa<0$ the sought surfaces become complex conjugate.

An alternative characterization of constant mean curvature surfaces in
the context of the isothermic transformation theory, via ``linear
conserved quantities'', has recently been given in
[\href[ref.bjrs08]{6, Sect~5}].  While allowing a uniform definition
of discrete constant mean curvature nets in spaces of constant
curvature and tying in very well with the transformation theory this
approach uses the full power of the isothermic transformations, hence
lacks the immediacy of the original ideas.

On the other hand, curvatures of circular surfaces with respect to
arbitrary Gauss maps based on Steiner's formula were introduced in
[\href[ref.sc03]{14, Sect 3}] and [\href[ref.sc06]{15}].  A curvature
theory for general polyhedral surfaces based on the notions of
parallel surfaces and mixed area is developed in
[\href[ref.bpw10]{5,~Def~8}], see also [\href[ref.bosu08]{4, Def
4.45}].  Embedding an ambient constant curvature geometry into a
linear space then allows to define these mixed area curvatures also
for discrete conjugate or principal nets in space forms, cf
[\href[ref.hrsy09]{12, Sect 6.3}], thereby leading to a simple
definition of discrete constant mean curvature surfaces in space
forms.  The necessity\footnote{The notion of mixed area relies on a
notion of parallel quadrilaterals or, more generally, polygons.} to
embed the ambient space form geometry into an affine space suggests
that, again, M\"obius geometry is a natural realm to work in.  A
purely sphere geometric treatment of constant mean curvature nets and,
more generally, discrete linear Weingarten surfaces can be obtained in
the realm of Lie sphere geometry, see [\href[ref.bjr11]{7}].

In the present paper we show that the discrete principal nets in a
quadric of constant curvature that have constant mixed area mean
curvature (see \href[thm.ma]{Def 2.3} and \href[def.maH]{3.2}) can be
characterized by the existence of a K\"onigs dual in a ``concentric
quadric'' (see \href[def.concentric]{Def 3.3}).  In this way, we
obtain a characterization that is very close to the original
definition of constant mean curvature surfaces in Euclidean space from
[\href[ref.bopi99]{2}] or [\href[ref.jhp99]{9}].  On the other hand,
we obtain the same class of ``constant mean curvature nets'' in space
forms as have been defined in [\href[ref.bjrs08]{6}] since they are
characterized by constant mixed area mean curvature, see
[\href[ref.bjr11]{7}].

Despite the fact that setup and some notions in the text are motivated
by M\"obius geometric considerations, the key arguments are of an
affine nature: both, a quadric of constant curvature (and its
concentric quadrics) as well as the K\"onigs dual, rely on a choice of
an affine subgeometry of the projective space that M\"obius geometry
is modelled on.  Thus no knowledge of M\"obius geometry will be
required to follow our main line of thought.  However, for full
appreciation of the geometric interrelations, some familiarity with
M\"obius geometry and, in particular, the theory of discrete
isothermic surfaces will be useful; for background material the
interested reader is referred to [\href[ref.imdg]{10}].

{\it Acknowledgements.\/} We would like to thank F Burstall and W
Rossman for fruitful and enjoyable discussions.  Thanks also go to the
referees whose comments greatly helped to improve the text.

\section K\"onigs duality.
Throughout the paper any occuring discrete net will be a (discrete)
conjugate net, that is, a quadrilateral net with planar faces.  For
simplicity we shall restrict to $\Z^2$ as a domain, though most of
what follows remains valid with only few adjustments when $\Z^2$ is
replaced by a more general quadgraph, that is, a $2$-dimensional cell
complex with quadrilateral $2$-cells.

Thus let $s:\Z^2\to\RP^4$ denote a discrete K\"onigs net, that is, a
discrete conjugate net so that, at each vertex $s_i$, $i\in\Z^2$, the
intersection points of pairs of diagonals of the four incident
quadrilaterals of $s_i$ are coplanar [\href[ref.bosu08]{4, Thm 2.26}],
see Fig 1.

Projecting to an affine subgeometry, the resulting K\"onigs net ${\fk
s}:\Z^2\to{\fk E}^4$ in affine space admits a dual K\"onigs net ${\fk
s^\ast}:\Z^2\to{\fk E}^4$:

\proclaim\htag[def.kdual]{2.1 Def}.
Two conjugate nets ${\fk s,s^\ast}:\Z^2\to{\fk E}^4$ in an affine
space ${\fk E}^4$ are called {\em K\"onigs dual} if they are
edge-parallel and have parallel non-corresponding diagonals.

In fact, the existence of such a dual ${\fk s^\ast}$ for a discrete
conjugate net ${\fk s}$ characterizes K\"onigs nets, cf
[\href[ref.bosu08]{4, Def 2.22}].

\vskip\abovedisplayskip
\hglue .035\hsize
\fig[36]Fig 1: K\"onigs net\hglue 4em -- 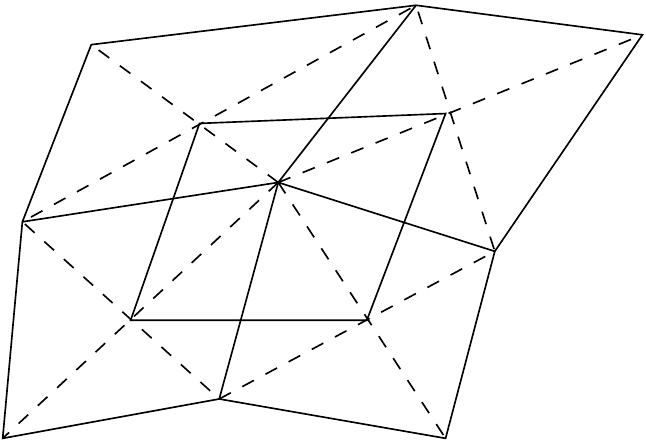
\put(0.0,0.0){}
\endfig\hfill
\fig[54]Fig 2: Dual quadrilaterals -- 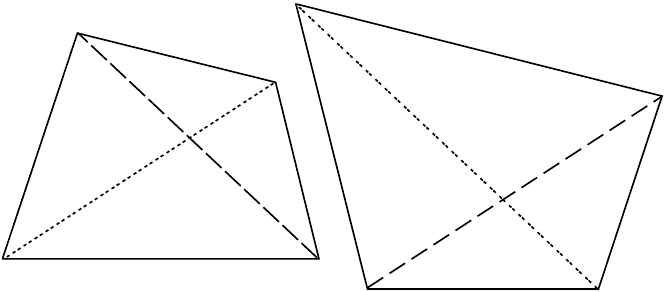
\put(0.00,0.03){${\fk s}_i$}
\put(0.45,0.03){${\fk s}_j$}
\put(0.38,0.36){${\fk s}_k$}
\put(0.11,0.42){${\fk s}_l$}
\put(0.50,0.01){${\fk s}^\ast_j$}
\put(0.91,0.01){${\fk s}^\ast_i$}
\put(0.99,0.33){${\fk s}^\ast_l$}
\put(0.40,0.44){${\fk s}^\ast_k$}
\endfig

In what follows we adopt the following notations, cf
[\href[ref.bosu08]{4, Def 2.23}]:

\proclaim\htag[def.dec]{2.2 Def}.
A {\em discrete $1$-form\/} is a map $(ij)\mapsto\omega_{ij}$ which
assigns values (in a vector space) to oriented edges, so that
$\omega_{ij}+\omega_{ji}=0$;\\ an {\em edge function\/} is a map
$(ij)\mapsto g_{ij}$ which assigns values to unoriented edges, that
is, $g_{ij}=g_{ji}$.

Given a map $i\mapsto g_i$ on vertices, its ``discrete differential''
is obtained by taking differences while the corresponding discrete
edge function is obtained by averaging: $$ (ij)\mapsto
dg_{ij}:=g_j-g_i \enspace{\rm and}\enspace (ij)\mapsto g_{ij}:={1\over
2}(g_i+g_j) $$ Note that $dg$ is a discrete $1$-form and that we
obtain a discrete Leibniz rule with these notations: $$
d(fg)_{ij}=df_{ij}\,g_{ij}+f_{ij}\,dg_{ij}.  $$ As diagonals of
elementary quadrilaterals $(ijkl)$ play a prominent role in our
analysis, it will be convenient to denote them by $ \delta g_{ik} :=
g_k-g_i.  $

\vskip\abovedisplayskip\hglue .35\hsize
\fig[30]Fig 3: Elementary quadrilateral -- 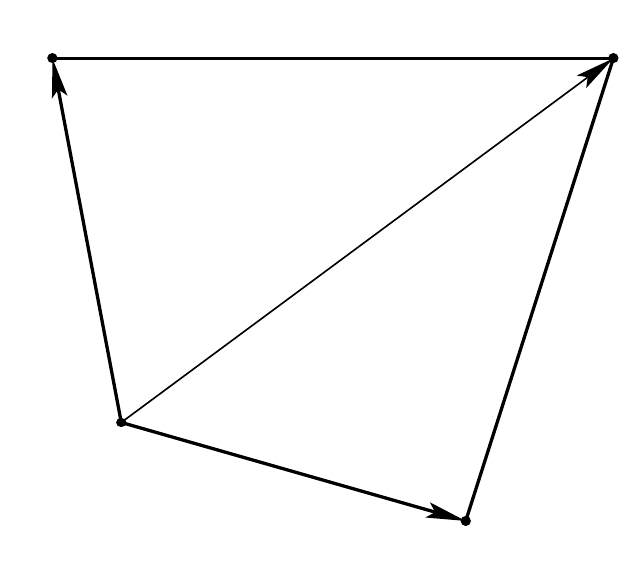
\put(0.15,0.20){${\fk s}_i$}
\put(0.75,0.05){${\fk s}_j$}
\put(0.99,0.85){${\fk s}_k$}
\put(0.00,0.85){${\fk s}_l$}
\put(0.45,0.22){$d{\fk s}_{ij}$}
\put(0.17,0.55){$d{\fk s}_{il}$}
\put(0.65,0.55){$\delta{\fk s}_{ik}$}
\endfig

Thus the K\"onigs dual ${\fk s}^\ast$ of a K\"onigs net ${\fk s}$ in
${\fk E}^4$ can be characterized by $$ d{\fk s}^\ast_{ij} \wedge d{\fk
s}_{ij} = 0 \enspace{\rm and}\enspace \delta{\fk s}^\ast_{ik} \wedge
\delta{\fk s}_{jl} = 0 $$ on all edges $(ij)$ and all elementary
quadrilaterals $(ijkl)$, respectively.  As $ \delta{\fk
s}_{ik}\wedge\delta{\fk s}^\ast_{jl} = \delta{\fk
s}^\ast_{ik}\wedge\delta{\fk s}_{jl} $ for edge-parallel nets ${\fk
s}$ and ${\fk s}^\ast$, the fact that non-corresponding diagonals of
${\fk s}$ and ${\fk s}^\ast$ are parallel can alternatively be
characterized by the vanishing of the mixed area\footnote{Note that
the presented notion of mixed area uses nothing but the linear
structure of ${\fk E}^4\cong\R^4$.} of corresponding elementary
quadrilaterals, cf [\href[ref.bosu08]{4, Thm 4.42}]:

\proclaim\htag[thm.ma]{2.3 Lemma \& Def}.  Two edge-parallel conjugate
nets ${\fk s,s^\ast}:\Z^2\to{\fk E}^4$ are K\"onigs dual if and only
if, on every elementary quadrilateral $(ijkl)$, their {\em mixed area}
$$ A({\fk s,s}^\ast)_{ijkl} : = {1\over4}\{ \delta{\fk
s}_{ik}\wedge\delta{\fk s}^\ast_{jl} + \delta{\fk
s}^\ast_{ik}\wedge\delta{\fk s}_{jl} \} = 0.  \eqno ma$$

To relate our current approach to discrete constant mean curvature
surfaces in space forms to the one given in [\href[ref.bjrs08]{6}] it
will be useful to describe the K\"onigs dual in terms of a Christoffel
formula, see [\href[ref.bosu08]{4, Thm 2.31}]: $$ d{\fk s}^\ast_{ij} =
{1\over\nu_i\nu_j}\,d{\fk s}_{ij}, \eqno christkoenigs$$ where the
``Christoffel symbol'' $\nu:\Z^2\to\R$ is a function that scales the
affine lift ${\fk s}$ of $s$ to a Moutard lift $ \mu = {{\fk
s}\over\nu}:\Z^2\to\R^5, $ that is, $\mu$ satisfies a discrete Moutard
equation\footnote{Note: while $s$ is a conjugate net in its projective
ambient geometry its Moutard lift $\mu$ may not have planar faces in
its linear ambient geometry.}, $\delta\mu_{ik}\wedge\delta\mu_{jl}=0$,
see [\href[ref.bosu08]{4, Thm 2.32}].

\section Constant mean curvature via K\"onigs duality.
Our aim is to classify K\"onigs nets in a quadric so that the K\"onigs
dual of an affine projection takes values in the same --- or, more
generally, a concentric --- quadric.  It turns out that this leads to
a characterization for constant mean curvature nets in space form.  A
natural problem to consider is a classification of the K\"onigs nets
that admit a K\"onigs dual in the same quadric: this captures those
constant mean curvature nets in space forms with a Lawson invariant
$H^2+\kappa>0$, in particular, all spherical constant mean curvature
nets and non-minimal constant mean curvature nets in Euclidean
geometry, cf [\href[ref.bjrs08]{6, Sect~4.3}].  To obtain a
classification of all constant mean curvature nets, including all
constant mean curvature nets in hyperbolic geometry, a wider class of
ambient quadrics for the K\"onigs dual needs to be considered: for
this reason we consider the ``concentric quadrics'' below.  A detailed
analysis of the occuring cases will be formulated in Sect~5 below.  As
a motivation we consider the case of non-minimal constant mean
curvature nets in Euclidean geometry:

{\bf Example}.  Let ${\fk Q}^3:=\{y\in\R^{4,1}\,|\,(y,y)=0,(y,q)=1\}$,
where ${\fk q}\in\R^{4,1}$ is isotropic, so that ${\fk Q}^3\cong\R^3$.
Namely, choose an ``origin'' ${\fk o}\in{\fk Q}^3$ (hence ${\fk o}$
and ${\fk q}$ span a Minkowski $\R^{1,1}\subset\R^{4,1}$) and observe
that $$ \R^3\cong\{{\fk o,q}\}^\perp\ni{\fk x}\mapsto {\fk s}={\fk
o+x}-{1\over2}({\fk x,x})\,{\fk q}\in{\fk Q}^3 \eqno stereoproj$$ is
an isometry (in the differential geometric sense), cf
[\href[ref.bosu08]{4, \S9.3.2}].  Suppose that ${\fk s}:\Z^2\to{\fk
Q}^3$ is an affine projection of a K\"onigs net in $\RP^4$ so that its
dual ${\fk s}^\ast:\Z^2\to{\fk Q}^3$ takes values in the same quadric
${\fk Q}^3$ of (vanishing) constant sectional curvature, $$ {\fk
o+x^\ast}-{1\over2}({\fk x^\ast,x^\ast})\,{\fk q} = {\fk s}^\ast
\enspace\Rightarrow\enspace d({\fk x^\ast-x})_{ij} \perp ({\fk
x^\ast-x})_{ij}, $$ using that $d{\fk
s}^\ast_{ij}={1\over\nu_i\nu_j}d{\fk s}_{ij}$.  Thus the stereographic
projections ${\fk x,x}^\ast:\Z^2\to\R^3$ of $s$ and $s^\ast$ are
K\"onigs dual circular nets\footnote{Therefore they are Christoffel
dual discrete isothermic nets (see below) in $\R^3$, cf
[\href[ref.bosu08]{4, Def 4.20}].} at a constant distance, that is,
they are parallel constant mean curvature nets in $\R^3$, with common
Gau\ss\ map ${\fk x-x}^\ast\over|{\fk x-x}^\ast|$, in the sense of
[\href[ref.jhp99]{9, Sect 5}] or [\href[ref.bosu08]{4, Cor 4.53}].
Thus requiring the dual of (an affine projection of) a K\"onigs net in
the quadric ${\fk Q}^3$ to take values in the same quadric forces (the
stereographic projection of) the discrete net to belong to the
subclass of constant mean curvature nets.  We shall see that this
breaking of symmetry is systemic.

Note that this construction yields a transformation\footnote{In the
case of discrete isothermic surfaces this is the ``Goursat
transformation'' of [\href[ref.imdg]{10, \S5.7.10}].} of discrete
K\"onigs nets in general quadrics: if the quadric is given as the
projective null-cone of a non-degenerate inner product on $\R^5$ then
a choice of a Minkowski plane, spanned by vectors ${\fk o,q}$ as
above, yields via ``stereographic projection'' a K\"onigs net ${\fk
x}:\Z^2\to\{{\fk o,q}\}^\perp$ whose dual in $\{{\fk o,q}\}^\perp$
lifts back to a K\"onigs net in the quadric.

We aim to provide a similar characterization as the one for constant
mean curvature nets in $\R^3$ above in the general case of an ambient
space form geometry.  Thus we consider the M\"obius quadric (see
[\href[ref.imdg]{10, Chap 1}] or [\href[ref.bosu08]{4, Sect 9.3}]) $$
S^3 = P({\cal L}^4), \enspace{\rm where}\enspace {\cal L}^4 =
\{y\in\R^{4,1}\,|\,(y,y)=0\} \eqno mobquad$$ and $(.,.)$ denotes a
Minkowski inner product on the space $\R^{4,1}\cong\R^5$ of
homogeneous coordinates of $\RP^4$.  A K\"onigs net $s:\Z^2\to
S^3\subset\RP^4$ then becomes circular: the intersection of
$S^3\subset\RP^4$ with the plane of each face yields the circumcircle
of that face.  That is, a discrete K\"onigs net in $S^3$ is a discrete
isothermic net, cf [\href[ref.bosu08]{4, Sect 4.3}]:

\proclaim\htag[def.isothermic]{3.1 Def}.
A discrete circular K\"onigs net is called {\em isothermic\/}.  A
discrete net in $S^3\subset\RP^4$ is {\em isothermic\/} if it is a
discrete K\"onigs net.

Alternatively, discrete isothermic nets can be characterized as
circular nets so that the (real) cross ratios $[{\fk s}_i,{\fk
s}_j,{\fk s}_k,{\fk s}_l]$ on elementary quadrilaterals $(ijkl)$ (as
cross ratios of four points on a conic) factorize, $$
[s_i,s_j,s_k,s_l] = {a_{ij}\over a_{jk}}, \eqno crisoth$$ where
$(ij)\mapsto a_{ij}$ is an edge labelling, that is, an edge function
($a_{ij}=a_{ji}$) that is constant across opposite edges
($a_{ij}=a_{kl}$), see [\href[ref.bosu08]{4, Thm 4.25}].

In consequence, the Christoffel formula \reqn{christkoenigs} can be
rewritten as (see [\href[ref.bosu08]{4, Thm 4.32}]) $$ d{\fk s}^\ast =
{a_{ij}\over|d{\fk s}_{ij}|^2}\,d{\fk s}_{ij}, \enspace{\rm
where}\enspace a_{ij}={|d{\fk s}_{ij}|^2\over\nu_i\nu_j}.  \eqno
christoffel$$

In the presence of a non-degenerate inner product, an affine
subgeometry is given by any choice of a non-zero ${\fk q}\in\R^{4,1}$
via $$ \R^4 \cong {\fk E}^4 : = \{{\fk y}\in\R^{4,1}\,|\,({\fk
y,q})=1\}, $$ and the projection $$ {\fk Q}^3 : = {\cal L}^4\cap{\fk
E}^4 \eqno cscquad$$ of $S^3$ to ${\fk E}^4$ is a quadric of constant
sectional curvature $\kappa=-({\fk q,q})$, cf [\href[ref.imdg]{10,
Sect 1.4}].

The characterization of K\"onigs duality in \href[thm.ma]{Lemma 2.3}
now suggests an intimate relation with the mean curvature of a
discrete isothermic net defined via mixed areas \reqn{ma}:

\proclaim\htag[def.maH]{3.2 Def}.
Let ${\fk s}:\Z^2\to{\fk Q}^3$ be a conjugate net; we say that a unit
vector field ${\fk n}:\Z^2\to S^{3,1}={\cal L}^4_1$ (see below) is a
{\em Gau\ss\ map\/} of ${\fk s}$ if ${\fk n}_i\in T_{{\fk s}_i}{\fk
Q}^3$ at all vertices $i\in\Z^2$ and ${\fk s}$ and ${\fk n}$ are edge
parallel, that is, if $$ ({\fk n,s})\equiv 0 \enspace{\rm and}\enspace
d{\fk n}_{ij} + k_{ij}d{\fk s}_{ij} = 0 \eqno rodrigues$$ on all edges
$(ij)$; the coefficients $k_{ij}$ yield the {\em edge principal
curvatures\/} of the pair $({\fk s,n})$.\\ Then the {\em Gau\ss} and
{\em mean curvature} $K$ and $H$, respectively, of the pair $({\fk
s,n})$ are defined by the equations $$ A({\fk n,n}) = K\,A({\fk s,s})
\enspace{\sl and}\enspace A({\fk n,s}) = -H\,A({\fk s,s}).  $$

A discrete conjugate net in a quadric ${\fk Q}^3$ of constant
curvature is circular, as discussed above, hence always admits a Gauss
map in the sense of this definition.  Note that both, mean and Gau\ss\
curvature, are defined on faces of a conjugate or isothermic net in
this way.

In the affine subgeometry of ${\fk E}^4$ it also makes sense to
consider ``concentric quadrics'', see Fig 4.

\proclaim\htag[def.concentric]{3.3 Def}.
Given an affine projection of $S^3$ to a quadric ${\fk Q}^3$ of
constant curvature in an affine subgeometry of $\RP^4$, we define its
{\em concentric quadrics} as $$ {\fk Q}^3_{r,t} := {\cal
L}^4_t\cap{\fk E}^4_r, \enspace{\sl where}\enspace {\fk E}^4_r :=
\{{\fk y}\in\R^{4,1}\,|\,({\fk y,q})=r\} \enspace{\sl and}\enspace
{\cal L}^4_t := \{y\in\R^{4,1}\,|\,(y,y)=t\}.  \eqno concentric$$

Note that we consider the full $2$-parameter family of quadrics ${\fk
Q}^3_{r,t}$ as ``concentric quadrics'', though we shall see below that
generically the family is essentially $1$-dimensional when using
suitable identifications.

These concentric quadrics lie in parallel affine subgeometries ${\fk
E}^4_r$, allowing their identification via projection along ${\fk q}$.
The terminology ``concentric'' is motivated by the case $({\fk
q,q})\neq0$, where the hyperplanes ${\fk E}^4_r$ inherit a
non-degenerate inner product and the points ${\fk c}_r:={r{\fk
q}\over({\fk q,q})}\in{\fk E}^4_r$ can be identified with a common
centre of the quadrics\footnote{In the positive definite case $({\fk
q,q})<0$ we must have $t({\fk q,q})<r^2$.  In the negative definite
case $({\fk q,q})>0$ we allow any values of $t$ and $r$, including the
case of the singular quadric with $t({\fk q,q})=r^2$.  In this way, we
obtain ``concentric'' quadrics of different types in the latter
case.}: $$ {\fk Q}^3_{r,t} = \{ {\fk y\in E}^4_r \,|\, ({\fk
y-c}_r,{\fk y-c}_r) = t-{r^2\over({\fk q,q})} =
{r^2+t\kappa\over\kappa} \}.  \eqno concentric1$$ Thus the quadrics
${\fk Q}^3_{r,t}$ and ${\fk Q}^3_{r',t'}$ with
$r^2+t\kappa=r'^2+t'\kappa$ can be identified via parallel projection
along ${\fk q}$.

In the degenerate case $({\fk q,q})=0$ the quadrics ${\fk Q}^3_{r,t}$
can simultaneously and isometrically be para\-metrized over the
Euclidean space $\R^3\cong\{{\fk o,q}\}^\perp$ by $$ {\fk Q}^3_{r,t} =
\{ r{\fk o + x} + {1\over 2r}(t-({\fk x,x})){\fk q}\,|\, {\fk
x}\in\R^3 \} \eqno concentric2$$ as long as $r\neq0$, and after
choosing\footnote{This choice of origin is geometrically irrelevant as
it only affects the parametrization up to translation.} an origin
${\fk o\in Q}^3={\fk Q}^3_{1,0}$ for the stereographic projection.
This shows that the ``concentric quadrics'' become a translation
family of paraboloids.  When $r=0$, on the other hand, ${\fk
Q}^3_{0,t}$ becomes a circular cylinder of radius $t>0$ and with
isotropic generators.

\vskip\abovedisplayskip
\hglue .13\hsize
\fig[73.87]Fig 4: Concentric quadrics -- 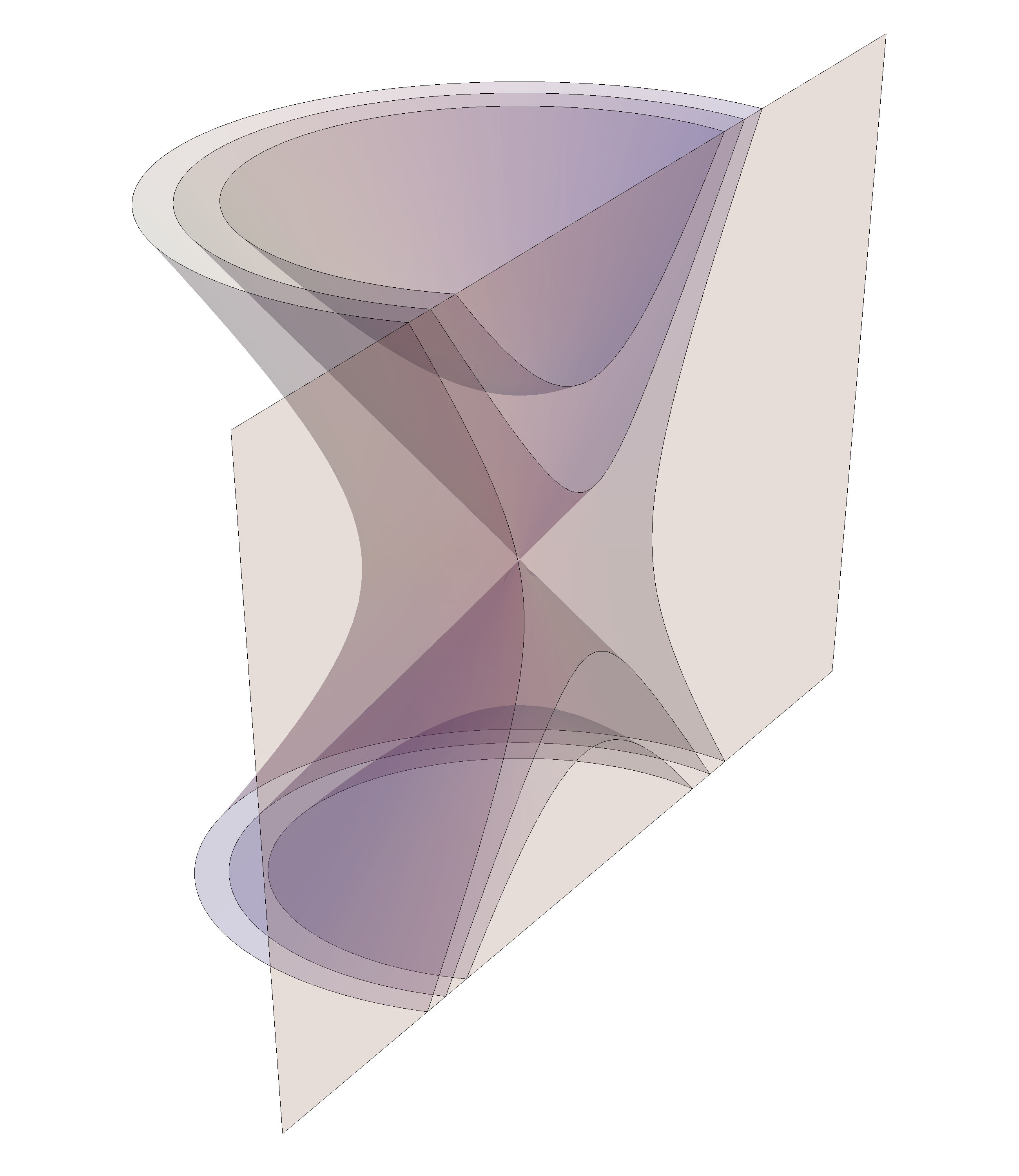
\put(0.83,0.52){${\fk E}_{1}^4$}
\put(0.12,1.00){${\cal L}_{1}^4$}
\put(0.18,1.00){${\cal L}_{0}^4$}
\put(0.24,1.00){${\cal L}_{-1}^4$}
\put(0.50,0.65){${\fk Q}_{1,1}^3$}
\put(0.50,0.75){${\fk Q}_{1,0}^3$}
\put(0.50,0.82){${\fk Q}_{1,-1}^3$}
\endfig

Now suppose that there is an affine projection ${\fk s}:\Z^2\to{\fk
Q}^3$ of our K\"onigs net $s:\Z^2\to S^3$, so that its K\"onigs
dual\footnote{In the cases ${\fk(q,q)}<0$, where the quadrics ${\fk
Q}^3_{r,t}$ are homothety equivalent, we could without loss of
generality restrict to ${\fk Q}^3$; in the cases ${\fk(q,q)}\geq0$
however, the quadrics come in different flavours so that a unified
analysis does require this more general setting.} $$ {\fk
s}^\ast:\Z^2\to{\fk Q}^3_{r,t} $$ takes values in one of the
concentric quadrics ${\fk Q}^3_{r,t}$.

As both $({\fk s,s}),({\fk s^\ast,s^\ast})\equiv \const$, we learn
that $ ({\fk s}_{ij},d{\fk s}_{ij}) = ({\fk s}^\ast_{ij},d{\fk
s}^\ast_{ij}) = 0 $ on any edge $(ij)$ of $\Z^2$; hence, $$ d({\fk
s,s}^\ast)_{ij} = (d{\fk s}_{ij},{\fk s}^\ast_{ij}) + ({\fk
s}_{ij},d{\fk s}^\ast_{ij}) = 0 $$ since ${\fk s}$ and ${\fk s}^\ast$
have parallel edges, which shows that:

\proclaim\htag[thm.distance]{3.4 Lemma}.
For K\"onigs dual nets ${\fk s}$ and ${\fk s}^\ast$ in concentric
quadrics we have $({\fk s,s}^\ast)\equiv \const$.\\ When $r=1$ then
${\fk s}$ and ${\fk s}^\ast$ lie in the same affine space ${\fk E}^4$
at constant distance $|{\fk s^\ast-s}|^2=t-2({\fk s,s^\ast})$.

In the positive definite case, where ${\fk s}$ and ${\fk s}^\ast$ take
values in two concentric spheres, this result can alternatively be
obtained by a pretty and more geometric argument: restricting focus to
the plane containing two corresponding edges $d{\fk s}_{ij}$ and
$d{\fk s}^\ast_{ij}$, we obtain a figure of two parallel segments in
concentric circles --- which therefore form an equilateral trapezoid.

\vskip\abovedisplayskip
\hglue .35\hsize
\fig[30]Fig 5: Constant distance in the sphere case
$|{\fk q}|^2<0$ -- 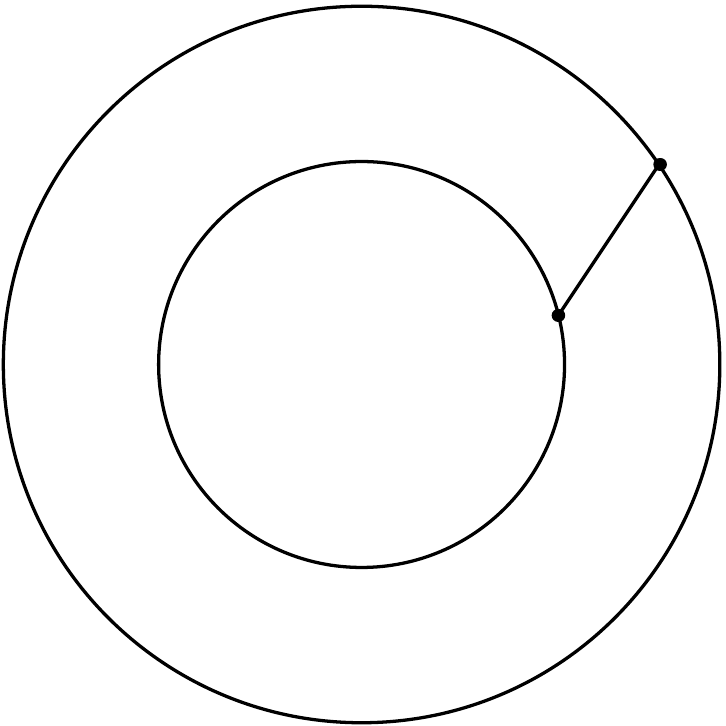
\put(0.94,0.05){${\fk E}_{1}^4$}
\put(-.15,0.50){${\fk Q}_{1,t}^3$}
\put(0.25,0.50){${\fk Q}^3={\fk Q}_{1,0}^3$}
\put(0.72,0.54){${\fk s}$}
\put(0.94,0.80){${\fk s}^\ast$}
\endfig

\proclaim\htag[thm.cmc]{3.5 Thm}.
If a K\"onigs net $s:\Z^2\to S^3\subset\RP^4$ admits an affine
projection ${\fk s}$ so that its K\"onigs dual ${\fk s}^\ast$ takes
values in a concentric\footnote{As the K\"onigs dual is translation
invariant ``concentric'' refers to the type of the target quadric
rather than its position.} quadric \reqn{concentric} then ${\fk s}$
has constant mean curvature $H$.

In order to prove Thm 3.5 we wish to show that the mixed area mean
curvature with respect to a suitable tangent plane congruence ${\fk
t}$ for ${\fk s}$ is constant.  Thus set $$ {\fk t} : = {\fk s}^\ast +
({\fk(s,s^\ast)(q,q)}-r)\,{\fk s} - {\fk(s,s^\ast)\,q} \eqno tpc$$ and
observe that ${\fk t\in\{s,q\}^\perp}=T_{{\fk s}}{\fk Q}^3$; in
particular, $({\fk t,t})>0$ unless ${\fk t}=0$ since ${\fk Q}^3$ is a
Riemannian quadric of constant curvature.  As $$ ({\fk t,t}) = t -
2r({\fk s,s^\ast}) + ({\fk q,q})({\fk s,s^\ast})^2 \equiv \const \eqno
tpcnorm$$ we learn that ${\fk t}_i\neq0$ at all points $i\in\Z^2$
except when ${\fk s}$ is ``self-dual'', that is, ${\fk s}^\ast$ is
obtained from ${\fk s}$ by a homothety and a translation.  Excluding
this self-dual case, ${\fk t}$ defines a suitable tangent plane
congruence of ${\fk s}$ in the sense that the pair $({\fk s,t})$
becomes a principal contact element net [\href[ref.bosu08]{4,
Def~3.24}] or discrete ``Legendre map'' [\href[ref.bjr11]{7, Def 3.1}]
as ${\fk t}$ and ${\fk s}$ have parallel edges: with the (edge)
principal curvatures $k_{ij}$ defined by $ 0 = d{\fk t}_{ij} +
k_{ij}\,d{\fk s}_{ij}, $ cf \reqn{rodrigues}, we obtain $$ \kappa_{ij}
: = {\fk t}_i + k_{ij}{\fk s}_i = {\fk t}_j + k_{ij}{\fk s}_j $$ as
the common (curvature) spheres associated to an edge.  Note that, in
the non-degenerate case $({\fk q,q})\neq0$ and denoting ${\fk
c}_r={r{\fk q}\over{\fk(q,q)}}$ as before, $$ {\fk t} = ({\fk
s^\ast-c}_r) + \{{\fk(s,s^\ast)(q,q)}-r\}({\fk s-c}_1).  \eqno
tpcnonflat$$

Thus the mixed area mean curvature of ${\fk s}:\Z^2\to{\fk Q}^3$ with
respect to ${\fk n}:={1\over\sqrt{\fk(t,t)}}\,{\fk t}$ as a Gau\ss\
map is now determined from $$ 0 = H\,A({\fk s,s}) + A({\fk s,n}) = \{
H + {{\fk(s,s^\ast)(q,q)}-r\over\sqrt{\fk(t,t)}} \}\, A({\fk s,s}) $$
since $A({\fk s,s^\ast})=A({\fk s,q})=0$.  Hence ${\fk s}$ has
constant mean curvature with respect to ${\fk n}$, proving
\href[thm.cmc]{Thm~3.5}.

\section K\"onigs duality via constant mean curvature.
To see that, conversely, every constant mean curvature net ${\fk s}$
in a quadric ${\fk Q}^3$ of constant sectional curvature, given in
\reqn{cscquad}, arises via K\"onigs duality in concentric quadrics we
shall construct a suitable dual for a given net ${\fk s}:\Z^2\to{\fk
Q}^3$ with Gau\ss\ map ${\fk n}$ that has constant mean curvature $H$
in the mixed area sense: $ 0 = A({\fk n,s}) + H\,A({\fk s,s}).  $

Thus we use the ``mean curvature sphere congruence''\footnote{We shall
convince ourselves below that ${\fk z}$ does indeed represent the mean
curvature sphere congruence of [\href[ref.bjrs08]{6, Def 5.1}].} of
the pair $({\fk s,n})$ to obtain a prototype K\"onigs dual for ${\fk
s}$:

\proclaim\htag[def.csc]{4.1 Lemma \& Def}.
Suppose the $({\fk s,n})$ has constant mean curvature $H$ in a quadric
${\fk Q}^3$ of constant curvature.  Then ${\fk s}:\Z^2\to{\fk E}^4$
has a K\"onigs dual ${\fk z}$ in the parallel affine subgeometry ${\fk
E}^4_H$ given by $$ {\fk z} : = {\fk n} + H{\fk s}:\Z^2\to{\fk
E}^4_H\subset\R^{4,1} \eqno csc$$ The sphere congruence $z=\R{\fk
z}:\Z^2\to\RP^4$ will be called the {\em mean curvature sphere
congruence} of $({\fk s,n})$.

Clearly $A({\fk z,s})=0$ and $({\fk z,q})=H$, so that ${\fk s}$ and
${\fk z}$ are K\"onigs dual nets in the parallel affine subgeometries
of ${\fk E}^4$ and ${\fk E}^4_H$, respectively, by \href[thm.ma]{Lemma
2.3}.  This proves Lemma 4.1.

The converse of \href[thm.cmc]{Thm 3.5} is now a straightforward
consequence:

\proclaim\htag[thm.csc]{4.2 Thm}.
Let ${\fk s}:\Z^2\to{\fk Q}^3$ be a discrete isothermic net in a
quadric of constant curvature $\kappa$ with constant mean curvature
$H$ with respect to a Gau\ss\ map ${\fk n}$.  Then its mean curvature
sphere congruence $z$ projects to a K\"onigs dual net ${\fk
s}^\ast={\fk z}$ in a concentric quadric ${\fk Q}^3_{r,t}$.

Namely, $({\fk z,z})=1$ so that ${\fk s}$ and ${\fk z}$ take values in
the concentric quadrics ${\fk Q}^3$ and ${\fk Q}^3_{H,1}$, proving the
theorem.

Note that the K\"onigs dual projection of the mean curvature sphere
congruence is not unique: as K\"onigs duality is invariant under
homotheties and translations, the nets $$ {\fk s}^\ast : = \lambda{\fk
z}+\mu{\fk q}:\Z^2\to{\fk Q}^3_{r,t}, \enspace{\rm where}\enspace
r=\lambda H-\mu\kappa \enspace{\rm and}\enspace
t=\lambda^2+2\lambda\mu H-\mu^2\kappa, \eqno centraldual$$ yield
K\"onigs duals into a family of concentric quadrics.  Identifying
parallel affine subgeometries via projection along ${\fk q}$ this is a
$1$-parameter family: in the non-flat cases $\kappa\neq0$ we learn
from \reqn{concentric1} that the translation part given by $\mu$ plays
no geometric role as $$ {\fk Q}^3_{r,t} = \{ {\fk y\in E}^4_r\,|\,
({\fk y-c}_r,{\fk y-c}_r) = \lambda^2{H^2+\kappa\over\kappa} \}, \eqno
concentricZ1$$ whereas this becomes clear from the parametrizations
given in \reqn{concentric2} in the flat case: $$ {\fk Q}^3_{r,t} = \{
\lambda H{\fk o+x} + {1\over 2H}(\lambda+2\mu H-({\fk x,x}))\,{\fk
q}\,|\, {\fk x}\in\R^3 \}.  \eqno concentricZ2$$ In the case
$H=\kappa=0$, where ${\fk s}$ is a discrete Euclidean minimal surface,
${\fk s}^\ast$ takes values in an isotropic cylinder of radius
$\lambda$ and axis ${\fk q}$ and projects to the Euclidean Gau\ss\ map
of the minimal net in $\R^3$.

Note that we obtain a K\"onigs dual net ${\fk s}^\ast$ in the original
quadric ${\fk Q}^3={\fk Q}^3_{1,0}$ from \reqn{centraldual} if and
only if the Lawson invariant $H^2+\kappa>0$, with $$ \lambda =
\pm{1\over\sqrt{H^2+\kappa}} \enspace{\rm and}\enspace \mu =
\mp{1\over\sqrt{H^2+\kappa}(H\pm\sqrt{H^2+\kappa})}.  $$

As we have now identified the K\"onigs nets in a quadric of constant
curvature that have a K\"onigs dual in a concentric quadric as the
mixed area constant mean curvature nets, it follows from
[\href[ref.bjr11]{7}] that these are precisely the constant mean
curvature nets in the sense of [\href[ref.bjrs08]{6}].

In fact, our ``mean curvature sphere congruence'' ${\fk z}$ is that of
[\href[ref.bjrs08]{6, Def 5.1}], cf [\href[ref.bjr11]{7, Expl 4.2}],
and in case of a positive Lawson invariant $H^2+\kappa>0$ the duals
${\fk s}^\ast:\Z^2\to{\fk Q}^3={\fk Q}^3_{1,0}$ obtained from
\reqn{centraldual} are the ``complementary nets'' of
[\href[ref.bjrs08]{6, Sect 4.3}]: special Darboux transforms of an
isothermic net\footnote{Recall, from [\href[ref.jhp99]{9}] or
[\href[ref.bosu08]{4, Def 4.27}], that a Darboux transform $\hat s$ of
an isothermic net $s:\Z^2\to S^3\subset\RP^4$ is defined by a cross
ratio condition on edges $(ij)$: $$ [s_i,s_j,\hat s_j,\hat s_i] =
\lambda a_{ij}, $$ where $\lambda\in\R$ is a parameter, and $a$ the
cross ratio factorizing edge labelling of \reqn{crisoth}.  In
particular, the Darboux transformation only depends on the conformal
ambient geometry of an isothermic net, and the defining cross ratio
condition is structurally the same as the characterization of
isothermic nets via factorizing cross ratios --- hence providing a
geometric realization of the $3d$-compatibility of the cross ratio
condition [\href[ref.bosu08]{4}].  A given isothermic net in $S^3$
admits a $4$-parameter family of Darboux transforms.} whose existence
is characteristic for discrete isothermic nets of constant mean
curvature, cf [\href[ref.jhp99]{9, Sect 5}] and our starting example.
We omit a detailed analysis to substantiate these claims here, as it
would require a discussion of the highly non-trivial techniques
introduced in [\href[ref.bjrs08]{6}], where the reader is referred for
further details.

However, as a consequence, examples and constructions of classes of
discrete constant mean curvature nets from the latter theory apply in
our approach: for example, the discrete minimal catenoid of
[\href[ref.jero10]{11}] or, more generally, the construction of
discrete constant mean curvature surfaces of revolution in space forms
provided there show that our theory is not empty.  A less involved
example of a dual pair of discrete constant mean curvature nets in the
$3$-sphere is discussed below.

\section Discrete cmc surfaces in space forms.
Thus we have obtained a characterization of discrete constant mean
curvature nets (in the sense of mixed areas) in space forms, that is,
simply connected complete spaces of constant sectional curvature (as
modeled by our quadrics ${\fk Q}^3={\fk Q}^3_{1,0}$), via K\"onigs
duality:

\proclaim\htag[thm.char]{5.1 Cor}.
A K\"onigs net $s:\Z^2\to S^3$ projects to a constant mean curvature
net ${\fk s}$ in a quadric ${\fk Q}^3$ of constant curvature if and
only if a K\"onigs dual of the projection takes values in a concentric
quadric ${\fk Q}^3_{r,t}$; the Gau\ss\ map ${\fk n}$ of ${\fk s}$ is
obtained by projecting the vector field joining ${\fk s}$ and its dual
${\fk s}^\ast$ onto the tangent bundle of ${\fk Q}^3$.

In the case $\kappa\neq0$ of a non-flat ambient geometry, the latter
assertion follows directly from \reqn{tpcnonflat} as $$ {\fk t} =
({\fk s^\ast-c}_r) - (r+\kappa({\fk s,s^\ast}))({\fk s-c}_1) $$ while,
in the case $\kappa=0$ of a flat ambient geometry, the parametrizations
\reqn{concentric2} reduce \reqn{tpc} to $$ {\fk t} = ({\fk
x}^\ast-r{\fk x}) - ({\fk x}^\ast-r{\fk x},{\fk x})\,{\fk q}.  $$

Constant mean curvatures surfaces in space forms come in different
types, depending on the sign of $H^2+\kappa$ or, equivalently,
depending on whether the Gau\ss\ equation reduces to the
$\sinh$-Gordon equation, the Liouville equation or the $\cosh$-Gordon
equation, respectively.  We shall see that these types are reflected
in the geometry of the K\"onigs dual ${\fk s}^\ast$ as well.  To see
this, we first note that $H^2+\kappa={r^2+t\kappa\over({\fk t,t})}$
using \reqn{tpcnorm}; in particular, the sign of $H^2+\kappa$ is the
same as that of $r^2+t\kappa$, showing that the ambient geometry of
the K\"onigs dual ${\fk s}^\ast$ is related to the sign of
$H^2+\kappa$, see \reqn{concentric1}.

\vskip\abovedisplayskip
\hglue .035\hsize
\fig[17]${\fk E}^4_1$, $\kappa>0$ -- 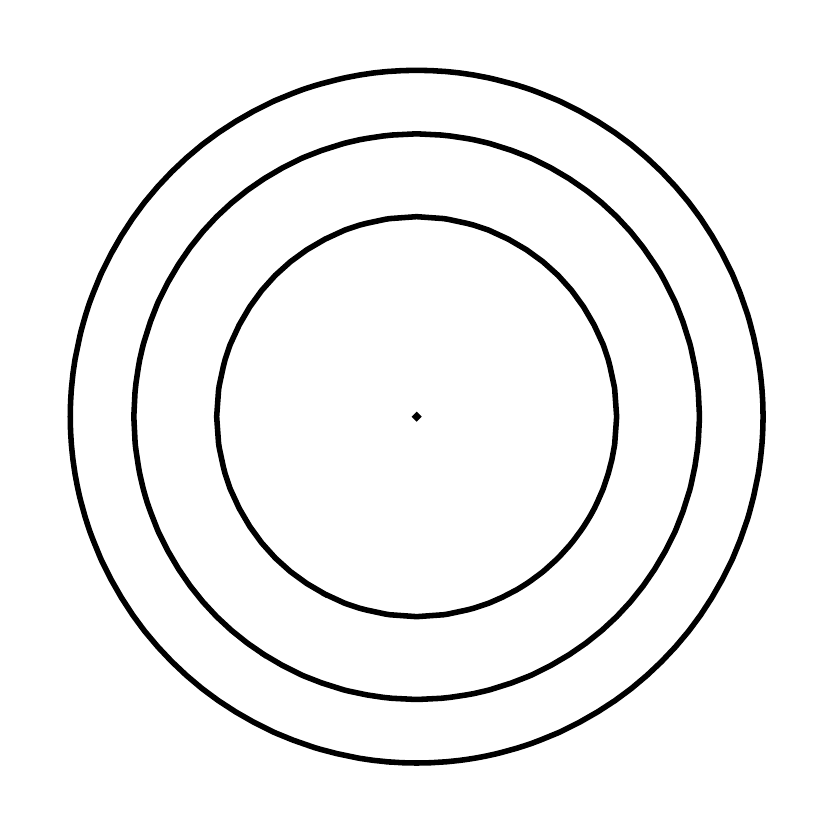
\endfig\hskip .02\hsize
\fig[17]${\fk E}^4_1$, $\kappa=0$ -- 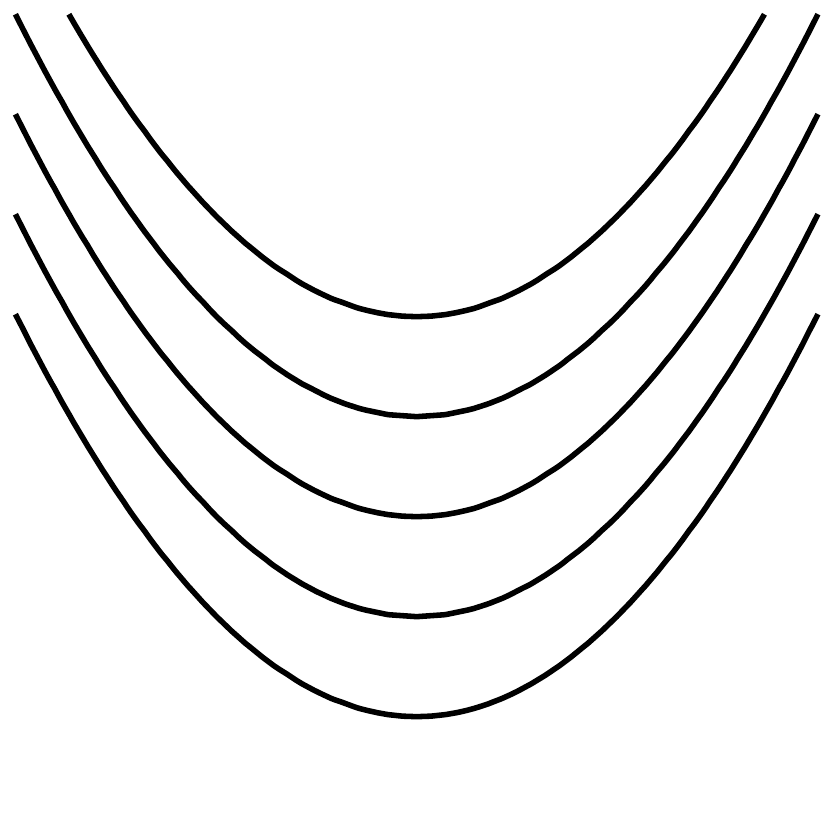
\endfig\hskip .02\hsize
\fig[17]${\fk E}^4_1$, $\kappa<0$ -- 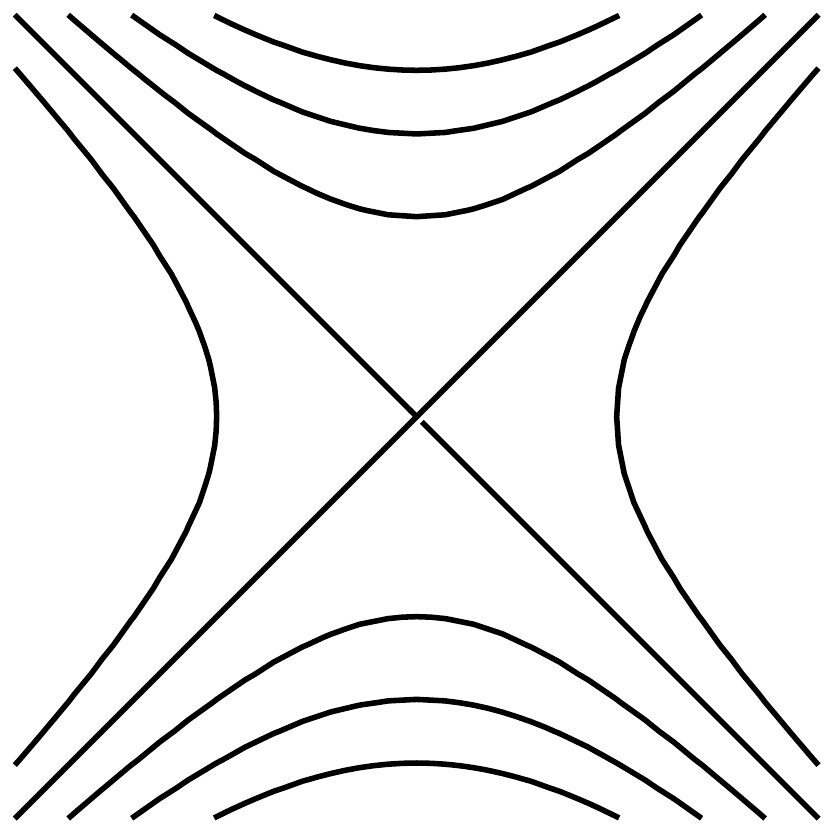
\endfig\hskip .02\hsize
\fig[17]${\fk E}^4_0$, $\kappa=0$ -- 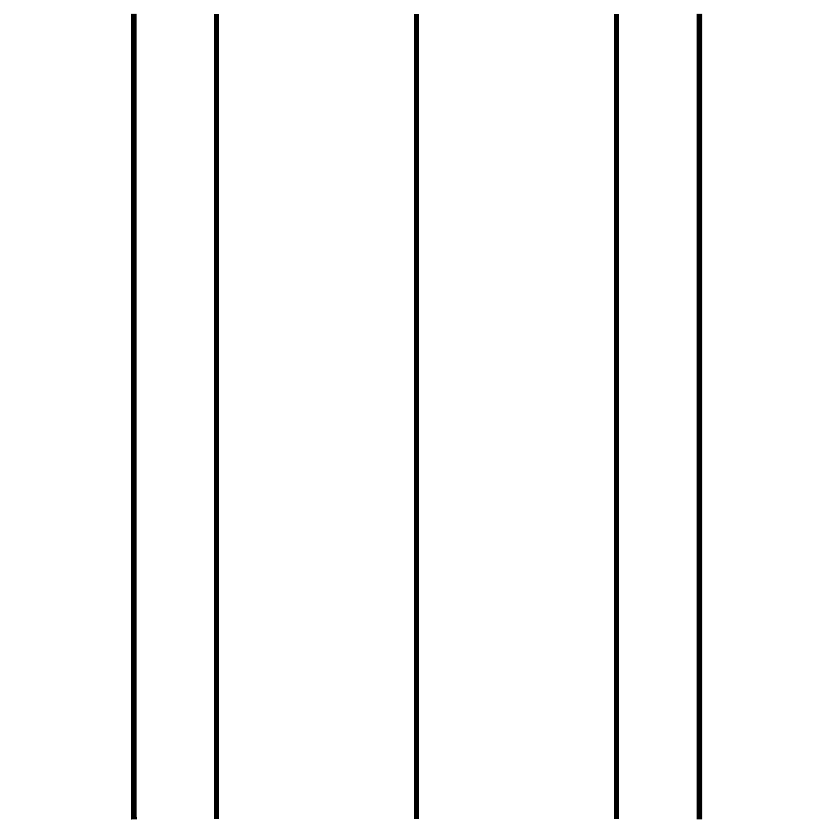
\endfig\hskip .02\hsize
\fig[17]${\fk E}^4_0$, $\kappa<0$ -- 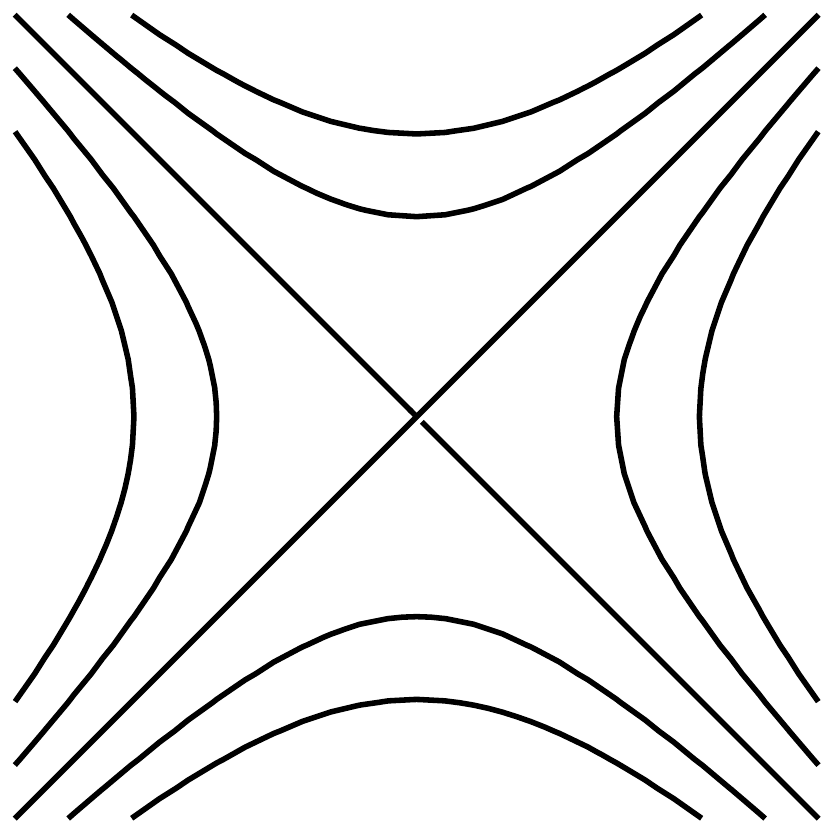
\endfig
\par\centerline{Fig 6: Concentric quadrics}

With this in mind we now provide a detailed analysis of the occuring
cases.

\item{(i)} $H^2+\kappa>0$: these surfaces arise in all space forms.

\item[4em]{(a)} Any constant mean curvature net in a spherical ambient
geometry, $\kappa>0$, belongs to this class as $H^2+\kappa>0$.  The
affine geometries of ${\fk E}^4_r$ are all Euclidean geometries and
\reqn{concentricZ1} shows that the quadrics ${\fk Q}^3_{r,t}$ are
concentric $3$-spheres.  Rescaling (in ${\fk E}^4_r$) and shifting
${\fk s}^\ast$ into ${\fk E}^4$, $$ {\fk s}^\ast \to
\pm{1\over\sqrt{r^2+t\kappa}}({\fk s^\ast-c}_r) + {\fk c}_1 =
\pm{1\over\sqrt{H^2+\kappa}}({\fk z}+{H\over\kappa}{\fk q}) -
{1\over\kappa}{\fk q}, \eqno s3shift$$ yields either of two antipodal
K\"onigs duals in ${\fk Q}^3$, that have constant (chordal) distances
$\sqrt{{2\over\kappa}(1\mp{H\over\sqrt{H^2+\kappa}})}$ from ${\fk s}$:
these form a pair of antipodal Darboux transforms of ${\fk s}$, whose
existence also characterizes discrete constant mean curvature nets in
a spherical geometry, see [\href[ref.bjrs08]{6, Sect 5.2}].

\item[4em]{(b)} In Euclidean geometry, $\kappa=0$, any non-minimal
constant mean curvature net belongs to the class under consideration.
The affine geometries of ${\fk E}^4_r$ inherit a degenerate metric
from $\R^{4,1}$ and from \reqn{concentricZ2} we learn that the
quadrics ${\fk Q}^3_{r,t}$ are a translation family of paraboloids.
Suitably rescaling and shifting ${\fk s}^\ast$ again, $$ {\fk s}^\ast
\to {1\over r}({\fk s}^\ast - {t\over 2r}{\fk q}) = {1\over H}({\fk
z}-{1\over 2H}{\fk q}), \eqno r3shift$$ places the K\"onigs dual in
${\fk Q}^3$ at a distance\footnote{In this case the chordal distance
coincides with the distance in $\R^3$.} ${1\over H}$ from ${\fk s}$:
using the parametrization $$ {\fk s} = {\fk o} + {\fk x} - {1\over
2}({\fk x,x})\,{\fk q} \enspace{\rm and}\enspace {\fk n} = {\fk g} -
({\fk g,x})\,{\fk q}, $$ of \reqn{concentric2}, where ${\fk g}$
denotes a Euclidean Gau\ss\ map of ${\fk x}$ in $\R^3=\{{\fk
o,q}\}^\perp$, we obtain $$ {\fk s}^\ast \simeq {\fk o} + ({\fk
x}+{1\over H}{\fk g}) - {1\over 2}(({\fk x,x}) +{2\over H}({\fk x,g})
+{1\over H^2})\,{\fk q}.  $$ Thus we recover the characterization of
constant mean curvature nets from [\href[ref.bopi99]{2, Sect 4}], or
from [\href[ref.jhp99]{9, Sect 5}] as this K\"onigs dual at constant
distance is also a Darboux transform of $s$.

\item[4em]{(c)} In the hyperbolic case, $\kappa<0$, the absolute value
of the mean curvature $H$ must be greater than that of a horosphere
for a constant mean curvature net to belong to the considered class.
Algebraically, the analysis of these nets is identical to the
spherical case, but their geometric interpretation becomes different:
now the affine geometries of ${\fk E}^4_{r,t}$ inherit a Minkowski
inner product from $\R^{4,1}$ and \reqn{concentricZ1} shows that the
quadrics ${\fk Q}^3_{r,t}$ are two-sheeted hyperboloids.  The
rescaling and shift from \reqn{s3shift} places the K\"onigs dual ${\fk
s}^\ast$ in either of the two sheets of ${\fk Q}^3$, so that they are
``antipodal'', that is, they can be identified via reflection in the
infinity boundary of ${\fk Q}^3$.  Consequently, the squares of their
constant chordal distances
${2\over\kappa}(1\mp{H\over\sqrt{H^2+\kappa}})$ to ${\fk s}$ have
different signs.  Again, both are Darboux transforms of $s$, placing
us in the situation analyzed in [\href[ref.bjrs08]{6, Sect 5.2}].

\item{(ii)} $H^2+\kappa=0$: for this to hold we must have
$\kappa\leq0$, that is, these surfaces only occur in Euclidean or
hyperbolic geometries.

\item[4em]{(b)} In the Euclidean case, $\kappa=0$, we consider minimal
nets in $\R^3$: the affine geometry of ${\fk E}^4_0$ inherits again a
degenerate inner product and the K\"onigs dual ${\fk s}^\ast$ of ${\fk
s}$ takes values in an isotropic cylinder ${\fk Q}^3_{0,t}\cong
S^2(\sqrt{t})\times\R{\fk q}$, $t>0$.  Thus, using the parametrization
from \reqn{concentric2} for ${\fk s}:\Z^2\to{\fk Q}^3={\fk
Q}^3_{1,0}$, $$ {\fk s} = {\fk o + x} - {1\over 2}({\fk x,x})\,{\fk q}
\enspace{\rm and}\enspace {1\over\sqrt{t}}({\fk s}^\ast-({\fk
s,s}^\ast)\,{\fk q}) = {\fk z} = {\fk g} - ({\fk g,x})\,{\fk q}, $$
where ${\fk g}:\Z^2\to S^2$ is a Euclidean Gau\ss\ map of ${\fk
x}:\Z^2\to\R^3\cong\{{\fk o,q}\}^\perp$ that renders ${\fk x}$ into a
minimal net in $\R^3$ in the sense of [\href[ref.bopi96]{1, Thm 8}].

\item[4em]{} Note that, while $|{\fk z-s}|^2=-2({\fk z,s})\equiv 0$,
shifting ${\fk s\to s-o}$ into the parallel affine subgeometry ${\fk
E}^4_0$ of ${\fk z}$, their ``distance'' $|{\fk z-s+o}|^2=(2{\fk
g-x,x})$ does not remain constant.  Also, in contrast to the case
$H^2+\kappa>0$, the K\"onigs dual ${\fk s}^\ast$ cannot be placed in
${\fk Q}^3$.  Consequently, it does not give rise to a Darboux
transform of $s$ in the sense of the isothermic transformation theory
of discrete isothermic nets in $S^3$.

\item[4em]{(c)} When $\kappa<0$ we are considering ``horospherical
nets'' in a hyperbolic geometry: ${\fk s}$ then has the same constant
mean curvature as a horosphere of the ambient hyperbolic geometry.
Now the affine subgeometry of ${\fk E}^4_r$ inherits a Minkowski inner
product from $\R^{4,1}$ and from \reqn{concentricZ1} we find that
${\fk Q}^3_{r,t}$ becomes the light cone of this Minkowski $4$-space.
Shifting both, ${\fk s}$ and ${\fk s}^\ast$ into the same affine
plane, say ${\fk E}^4_0$: $$ {\fk s\to s-c}_1 \enspace{\rm
and}\enspace {\fk s^\ast\to s^\ast-c}_r, $$ their distance $ |({\fk
s^\ast-c}_r)-({\fk s-c}_1)|^2 = {1-2r\over\kappa}-2({\fk s,s^\ast}) $
in ${\fk E}^4_0$ is constant by \href[thm.distance]{Lemma 3.4}.  On
the other hand, ${\fk s},{\fk s-c}_r:\Z^2\to{\cal L}^4$ yield maps
into $S^3\subset\RP^4$, allowing to interpret them as Darboux
transforms of each other --- which recovers the characterization of
horospherical nets of [\href[ref.imdg]{10, \S5.7.37}] via the
existence of a ``hyperbolic Gau\ss\ map'' as a totally umbilic Darboux
transform.

\item{(iii)} $H^2+\kappa<0$: these surfaces only occur in hyperbolic
geometry and include the interesting case of hyperbolic minimal
surfaces, such as the minimal catenoid of [\href[ref.jero10]{11}].

\item[4em]{(c)} As $\kappa<0$, the affine ambient geometry of the
K\"onigs dual ${\fk s}^\ast$ of ${\fk s}$ inherits again a Minkowski
inner product.  The target quadric ${\fk Q}^3_{r,t}$ of ${\fk s}^\ast$
now becomes a one-sheeted hyperboloid, see \reqn{concentricZ1}.
Arbitrarily rescaling ${\fk s}^\ast$ (in ${\fk E}^4_r$) and shifting
it into the affine ambient space ${\fk E}^4$ of ${\fk s}$, $$ {\fk
s}^\ast \to \lambda({\fk s-c}_r) + {\fk c}_1, $$ we obtain a K\"onigs
dual pair in ${\fk E}^4$ at a constant distance.  However, as in the
case of minimal nets in Euclidean geometry, ${\fk s}^\ast$ cannot be
placed in (an affine projection of) the M\"obius quadric
$S^3\subset\RP^4$, so that an interpretation of the dual as a Darboux
transformation of $s:\Z^2\to S^3$ is not available.

Thus, while in Euclidean and hyperbolic ambient geometries constant
mean curvature nets occur in different flavours, a beautiful and
simple result emerges in spherical ambient geometry:

\proclaim\htag[thm.spherical]{5.2 Cor}.
A K\"onigs net ${\fk s}:\Z^2\to S^3$ possesses a unit normal field
${\fk n}:\Z^2\to S^3$ so that $({\fk s,n})$ has constant mean
curvature if and only if ${\fk s}$ admits a (pair of antipodal)
K\"onigs dual(s) in $S^3\subset\R^4$.

{\bf Example}.  A simple example, beautifully illustrating this
corollary, is given by the Clifford tori in $S^3$.  Fixing a timelike
vector ${\fk q}\in\R^{4,1}$, $({\fk q,q})=-1$, the affine hyperplane
${\fk E}^4={\fk E}^4_1\cong\R^4$ and ${\fk Q}^3={\fk Q}^3_{1,0}\cong
S^3$ becomes the unit sphere in this Euclidean $4$-space.  Now let
$\varphi,\psi:\Z\to\R$ be two (strictly increasing) functions that
project to periodic $\R/(2\pi\Z)$-valued functions, that is,
$\varphi(m+m_0)=\varphi(m)+2\pi$ and $\psi(n+n_0)=\psi(n)+2\pi$ for
some $m_0,n_0\in\N$, let $\alpha\in(0,{\pi\over2})$, and
define\footnote{Note that, for simplicity, we omit the shift of ${\fk
s}$ into the affine plane ${\fk E}^4={\fk E}^4_1$ so that $|{\fk
s}|^2=1$ instead of taking values in the light cone of $\R^{4,1}$.  As
${\fk n}$ does take values in the linear space $\R^4={\fk
q}^\perp\subset\R^{4,1}$, i.e., is not shifted relative to the
$\R^{4,1}$-picture, it is a Gau\ss\ map of ${\fk s-q}:\Z^2\to{\fk
Q}^3$.} $$ {\fk s} : = \left({e^{i\varphi}\cos\alpha\atop
e^{i\psi}\sin\alpha}\right) \enspace{\rm and}\enspace {\fk n} : =
\left({-e^{i\varphi}\sin\alpha\atop
\phantom{-}e^{i\psi}\cos\alpha}\right), $$ where we identify
$\R^4\cong\C^2$.  Clearly, ${\fk s,n}:\Z^2\to S^3$ and ${\fk
n}\perp{\fk s}$.  Moreover, ${\fk s}$ and ${\fk n}$ are edge parallel,
$$ 0 = (d{\fk n} + \tan\alpha\,d{\fk s})_{(m,n)(m+1,n)} = (d{\fk n} -
\cot\alpha\,d{\fk s})_{(m,n)(m,n+1)}, $$ when $\varphi=\varphi(m)$ and
$\psi=\psi(n)$, so that ${\fk s}$ has (edge) principal curvatures
$\tan\alpha$ and $-\cot\alpha$ with respect to ${\fk n}$ as a Gau\ss\
map.  To determine the (face) mean curvature note that ${\fk s}$ and
$$ {\fk z} : = {\fk n}-\cot2\alpha\,{\fk s} =
{1\over\sin2\alpha}\left( -e^{i\varphi}\cos\alpha\atop
\phantom{-}e^{i\psi}\sin\alpha\right) $$ have parallel opposite
diagonals, $$ \delta{\fk z}_{(m,n)(m+1,n+1)} \parallel \delta{\fk
s}_{(m+1,n)(m,n+1)} \enspace{\rm and}\enspace \delta{\fk
z}_{(m+1,n)(m,n+1)} \parallel \delta{\fk s}_{(m,n)(m+1,n+1)}, $$
showing that $$ 0 = A({\fk z,s}) = A({\fk n,s}) - \cot2\alpha\,A({\fk
s,s}), \enspace{\rm hence}\enspace H \equiv -\cot2\alpha.  $$ Note
that, in this case, we obtain the mean curvature as the arithmetic
mean of the principal curvatures of the edges of a face --- in
general, this relation is far more intricate, cf [\href[ref.bpw10]{5,
Prop 11}].  Further, two (antipodal) dual nets of ${\fk s}$ in
$S^3\subset\R^4\cong{\fk E}^4$ are obtained as $$ {\fk s}^\ast =
\pm\sin2\alpha\,{\fk z} = \pm(\sin2\alpha\,{\fk n} - \cos2\alpha\,{\fk
s}).  $$ Note that in the case $\alpha={\pi\over 4}$ of a minimal
(square) Clifford torus ${\fk s}^\ast=\pm{\fk n}=\pm{\fk z}$.  If
$\varphi$ or $\psi$ projects to a periodic $\R/(\pi\Z)$-valued
function, with periods $m_1={m_0\over 2}$ or $n_1={n_0\over 2}$, then
the minimal Clifford torus becomes ``self-dual'' in the sense that $$
{\fk s}^\ast_{m,n} = \pm{\fk s}_{m+m_1,n} \enspace{\rm or}\enspace
{\fk s}^\ast_{m,n} = \mp{\fk s}_{m,n+n_1}.  $$

\goodbreak\vskip 3ex{\useFnt[cmbx10 scaled 1200]References}\vglue 1ex
{\frenchspacing
\item{\htag[ref.bopi96]{1.}} A Bobenko, U Pinkall: {\it Discrete
isothermic surfaces\/}; J reine angew Math 475, 187--208 (1996)
\item{\htag[ref.bopi99]{2.}} A Bobenko, U Pinkall: {\it Discretization
of surfaces and integrable systems\/}; Oxf Lect Ser Math Appl 16,
3--58 (1999)
\item{\htag[ref.bhs06]{3.}} A Bobenko, T Hoffmann, B Springborn: {\it
Minimal surfaces from circle patterns: geometry from combinatorics\/};
Ann Math 164, 231--264 (2006)
\item{\htag[ref.bosu08]{4.}} A Bobenko, Y Suris: {\it Discrete
differential geometry.  Integrable structure\/}; Grad Stud Math 98,
Amer Math Soc, Providence RI (2008)
\item{\htag[ref.bpw10]{5.}} A Bobenko, H Pottmann, J Wallner: {\it A
curvature theory for discrete surfaces based on mesh parallelity\/};
Math Ann 348, 1--24 (2010)
\item{\htag[ref.bjrs08]{6.}} F Burstall, U Hertrich-Jeromin, W
Rossman, S Santos: {\it Discrete surfaces of constant mean
curvature\/}; RIMS Kyokuroku Bessatsu 1880, 133--179 (2014)
\item{\htag[ref.bjr11]{7.}} F Burstall, U Hertrich-Jeromin, W Rossman:
{\it Discrete linear Weingarten surfaces\/}; arXiv:1406.1293 (2014)
\item{\htag[ref.gbpo96]{8.}} K Gro\ss{}e-Brauckmann, K Polthier: {\it
Numerical examples of compact constant mean curvature surfaces\/}; in
B Chow et al: Elliptic and parabolic methods in geometry 23--46 (1996)
\item{\htag[ref.jhp99]{9.}} U Hertrich-Jeromin, T Hoffmann, U Pinkall:
{\it A discrete version of the Darboux transform for isothermic
surfaces\/}; Oxf Lect Ser Math Appl 16, 59--81 (1999)
\item{\htag[ref.imdg]{10.}} U Hertrich-Jeromin: {\it Introduction to
M\"obius differential geometry\/}; London Math Soc Lect Note Series
300, Cambridge Univ Press, Cambridge (2003)
\item{\htag[ref.jero10]{11.}} U Hertrich-Jeromin, W Rossman: {\it
Discrete minimal catenoid in hyperbolic $3$-space\/}; Electr Geom
Model 2010.11.001 (2010)
\item{\htag[ref.hrsy09]{12.}} T Hoffmann, W Rossman, T Sasaki, M
Yoshida: {\it Discrete flat surfaces and linear Weingarten surfaces in
hyperbolic $3$-space\/}; Trans AMS 364, 5605--5644 (2012)
\item{\htag[ref.pipo93]{13.}} U Pinkall, K Polthier: {\it Computing
discrete minimal surfaces and their conjugates\/}; Experimental Math
2, 15--36 (1993)
\item{\htag[ref.sc03]{14.}} W Schief: {\it On the unification of
classical and novel integrable surfaces II.  Difference geometry\/};
Proc Royal Soc London A 459, 2449--2462 (2003)
\item{\htag[ref.sc06]{15.}} W Schief: {\it On a maximum principle for
minimal surfaces and their integrable discrete counterparts\/}; J Geom
Phys 56, 1484--1495 (2006)
\par}\vfill
\bgroup
\small
\def\addwd{\hsize=.35\hsize}
\def\udo{\vtop{\addwd
U Hertrich-Jeromin\\
Technische Universit\"at Wien\\
Wiedner Hauptstra\ss{}e 8--10/104\\
1040 Wien (Austria)\\
Email: udo.hertrich-jeromin@tuwien.ac.at
}}
\def\sasha{\vtop{\addwd
A Bobenko \& I Lukyanenko\\
Department of Mathematics\\
Technische Universit\"at Berlin\\
10623 Berlin (Germany)\\
Email: bobenko@math.tu-berlin.de\\
\phantom{Email: }lukyanen@math.tu-berlin.de
}}
\hbox to \hsize{\hfil \sasha \hfil \udo \hfil}
\egroup
\end{document}